# A GENERAL NONCONVEX LARGE DEVIATION RESULT II

By A. de Acosta

*Case Western Reserve University*

We refine the conditions for the lower bound in an abstract large deviation result with nonconvex rate function we had previously introduced. We apply the results to certain stochastic recursive schemes.

**1. Introduction.** In the recent paper [5], we introduced an abstract scheme designed to handle a broad class of large deviation problems in which the random variables take values in a topological vector space $E$ and the rate function is not convex. A rough description of our scheme is as follows. Let $E$ be as above, let $E^*$ be its dual space and let $\{Y_n\}_{n \in N}$ be $E$-valued random vectors. Assume:

(i) For certain functions $\Phi_n : E \times E^* \to \mathbf{R}$, all $n \in N$, all $\xi \in E^*$,

$$\mathbf{E} \exp[\langle Y_n, \xi \rangle - \Phi_n(Y_n, \xi)] = 1.$$

(ii) For a certain function $\Phi : E \times E^* \to \mathbf{R}$, all $x \in E$, all $\xi \in E^*$,

$$\lim_n n^{-1} \Phi_n(x, n\xi) = \Phi(x, \xi).$$

(iii) $\{\mathcal{L}(Y_n)\}_{n \in N}$ is exponentially tight.

Then under suitable regularity conditions on $\Phi$, $\{\mathcal{L}(Y_n)\}_{n \in N}$ satisfies the large deviation principle with rate function $\Phi^*(x, x)$, where for $x, y \in E$,

$$\Phi^*(x, y) = \sup_{\xi \in E^*} [\langle y, \xi \rangle - \Phi(x, \xi)].$$

Precise conditions under which the scheme is valid are given in Theorems 2.1 (upper bound) and 2.2 (lower bound) of [5].

While most conditions in Theorems 2.1 and 2.2 of [5] are formulated directly in terms of $\Phi$ and appear to be reasonably simple to verify in applications, condition (11) of Theorem 2.2 of [5]—an assumption on $\Phi^*$ involving sub-differentials—is in general more difficult to check (see the beginning of the proof of Theorem 2.2 below for a detailed statement of this condition). As is well known, a change of measure in some form is crucial in many proofs of large deviation lower bounds. The purpose of the condition is to ensure that for every point $x$ in the domain of the rate function, there exists a nearby "smooth point" $y$ such that the function values are close and there exists a suitable change of measure centered at $y$. If $E$ is a Banach space and $\Phi$ does not depend on $x$, the abundance of "smooth points" can be obtained from the Brondsted–Rockafellar theorem (see, e.g., [1]) and has been used, for example, in [4] in the proof of lower bounds with the convex rate function $\Phi^*$. Condition (11) of Theorem 2.2 in [5] may be regarded as a "nonconvex" version of the conclusion of the Brondsted–Rockafellar theorem.

The main objective of the present paper is to refine the abstract lower bound result in [5] by providing sufficient conditions for the subdifferentiability assumption in [5] which do not involve subdifferentials and are substantially easier to verify. What we prove in this context may be regarded as a "nonconvex" version of the Brondsted–Rockafellar theorem, guaranteeing the existence of an abundance of "smooth points." The tools used in the proof are a result of Zabell [17] on Mosco convergence of convex functions in locally convex spaces, the Schauder–Tychonoff fixed point theorem (see, e.g., [14]) and our recent result on dominating points of convex open sets in the context of general convex functions [6]. This objective is pursued in Section 2, Theorem 2.2. In Theorem 2.1 we present a simple improvement of the upper bound result in [5] in the framework of the present paper.

In Section 3 we present some applications of Theorems 2.1 and 2.2. In Theorem 3.1 we give a new approach to the study of large deviations for a recursive scheme based on an i.i.d. sequence of random vector fields, a question considered by Dupuis and Ellis [7] [see part 1 of Remark 3.2 ]. In Theorem 3.11 we consider the case when the recursive scheme is a stochastic Euler-type polygonal scheme for a dynamical system (see Remark 3.13).

We close this introduction with some remarks about the connection of our abstract scheme to certain items in the literature. A number of results on large deviations for trajectories of Markov processes [7, 9, 16, 15] involve nonconvex rate functions which are in fact of the form $\Phi^*(x,x)$ described above, although this aspect is not mentioned and $\Phi$ is not introduced. We showed in [5], Theorem 3.1, how our scheme applies to large deviations for the trajectories of a broad class of Markov processes, with the rate function initially given in the form $\Phi^*(x,x)$ and subsequently identified in a classical integral form. (Incidentally, the proof of that theorem can be simplified using the results of the present paper.) Our scheme is also related to ideas developed for the study of large deviations for semimartingales in [11] and [13], which present a general framework for the problem. However, the full details of the technical connection between this development and our scheme have yet to be elucidated; this remark applies as well to the Markov case mentioned above.

**2. The general large deviation results.** Throughout the section we assume

$E$ is a Banach space, $\mathcal{E}$ is the $\sigma$-algebra generated by the balls, $F$ is a subspace of $E^*$ such that $\langle \cdot, \xi \rangle$ is $\mathcal{E}$-measurable for every $\xi \in F$.

We have adopted this framework, which is less general than that in [5], in order to maintain some consistency in the presentation; in fact, it is only the proof of the lower bound that requires it. Two important cases are covered by these assumptions:

1. $E = C([0,1], \mathbf{R}^d), \mathcal{E} =$ Borel $\sigma$-algebra, $F = M([0,1], \mathbf{R}^d)$, the space of finite $\mathbf{R}^d$-valued vector measures on $[0,1]$.
2. $E = D([0,1], \mathbf{R}^d)$ endowed with the uniform norm, $\mathcal{E} = \sigma$-algebra generated by the evaluation maps, $F = M([0,1], \mathbf{R}^d)$.

The setting 2 was used in the application to stochastic equations in [5], and previously in [4].

For a function $\Phi : E \times F \to \mathbf{R}$, we define, for $x, y \in E$,

$$\Phi^*(x,y) = \sup_{\xi \in F} [\langle y, \xi \rangle - \Phi(x, \xi)].$$

In what follows, $\{a_n\}_{n \in N}$ is a positive sequence with $\lim_n a_n = \infty$.

THEOREM 2.1. *Let $\Phi_n, \Phi : E \times F \to \mathbf{R}$ be such that:*

1. *For all $\xi \in F, \Phi_n(\cdot, \xi)$ is $\mathcal{E}$-measurable.*
2. *For all $\xi \in F, \Phi(\cdot, \xi)$ is $\mathcal{E}$-measurable, continuous and satisfies $\Phi(x, 0) = 0$ for all $x \in E$.*
3. *For all $\xi \in F$, all compact sets $K \subset E$,*

$$b_n(K, \xi) \triangleq \sup_{x \in K} |a_n^{-1} \Phi_n(x, a_n \xi) - \Phi(x, \xi)| \to 0 \qquad \text{as } n \to \infty.$$

*For each $n \in \mathbf{N}$, let $Y_n$ be an $E$-valued, $\mathcal{E}$-random vector defined on $(\Omega_n, \mathcal{A}_n, \mathbf{P}_n)$ and assume:*

4. *For all $n \in \mathbf{N}, \xi \in F$,*

$$\mathbf{E}_n \exp[\langle Y_n, \xi \rangle - \Phi_n(Y_n, \xi)] = 1.$$

5. *$\{\mathcal{L}_{P_n}(Y_n)\}_{n \in N}$ is exponentially tight.*

*Then if* 1–4 *are satisfied, for every compact set* $K \subset E$,

$$\limsup_n a_n^{-1} \log \mathbf{P}\{Y_n \in K\} \leq - \inf_{x \in K} \Phi^*(x, x),$$

*and if* 1–5 *are satisfied, for every* $A \in \mathcal{E}$,

$$\limsup_n a_n^{-1} \log \mathbf{P}\{Y_n \in A\} \leq - \inf_{x \in \bar{A}} \Phi^*(x, x).$$

We omit the proof, which involves an easy modification of the proof of Theorem 2.1 of [5]. Note that condition 3 improves the corresponding assumption in Theorem 2.1 of [5]: for all $\xi \in F$,

$$\sup_{x \in E} |a_n^{-1} \Phi_n(x, a_n \xi) - \Phi(x, \xi)| \to 0 \qquad \text{as } n \to \infty.$$

Moreover, compared to Theorem 2.1 of [5], we are taking here $Z_n = Y_n$ and assumption 7 there is unnecessary.

For the main result, the large deviation lower bound, we further specify the framework as follows:

$$F = E_0^* \qquad \text{where } E_0 \text{ is a closed separable subspace of } E.$$

We need this assumption to ensure the applicability of the result on Mosco convergence in [17]. Note that the cases 1 and 2 mentioned above are still covered; in case 2, we take $E_0 = C([0,1], \mathbf{R}^d)$.

Recall that a function $\phi: F \to \mathbf{R}$ is $E_0$-*Gâteaux differentiable* at $\xi \in F$ if there exists a point $\triangledown \phi(\xi) \in E_0$ such that, for all $\eta \in F$,

$$\langle \triangledown \phi(\xi), \eta \rangle = \lim_{t \to 0} t^{-1} [\phi(\xi + t\eta) - \phi(\xi)].$$

Throughout the paper, when the gradient operator is applied to a function of two variables, it will refer to differentiation with respect to the second variable.

We use the notation $\partial \Phi^*(x, y)$ for the subdifferential of the convex function $\Phi^*(x, \cdot)$ at $y \in E$ (for the definition of subdifferential, see, e.g., [8]).

For $g: E \to \mathbf{R}^+, a \geq 0$, let $L(g, a) = \{x \in E : g(x) \leq a\}$.

THEOREM 2.2. *Assume that the hypotheses of Theorem* 2.1 *hold, and furthermore:*

6. *For all* $\xi \in F$,

$$\limsup_n \sup_{x \in E} |a_n^{-1} \Phi_n(x, a_n \xi)| < \infty.$$

7. *For every* $x \in E$ *such that* $\Phi^*(x, x) < \infty$, *there exists a neighborhood* $U$ *of* $x$ *such that, for all* $a \geq 0$,

$$\bigcup_{y \in U} L(\Phi^*(y, \cdot), a) \text{ is a relatively compact subset of } E_0.$$

8. If $x_n \to x$ in $E_0$ and $\xi_n \xrightarrow[w^*]{} \xi$ in $F$, then
$$\Phi(x,\xi) \leq \liminf_n \Phi(x_n, \xi_n).$$

9. For all $x \in E, \Phi(x, \cdot)$ is convex and $E_0$-Gâteaux differentiable on $F$. Moreover, for all $\xi \in F, x \in E_0, \phi'_{x,\xi}$ is continuous, where for $t \in \mathbf{R}$,
$$\phi_{x,\xi}(t) \triangleq \Phi(x, t\xi).$$

10. For all $\xi \in F$, the equation $x = \nabla \Phi(x, \xi)$ has at most one solution in $E_0$.
11. For all $x_0$ such that $\Phi^*(x_0, x_0) < \infty$, for every $\varepsilon > 0$, there exists $y_0 \in B(x_0, \varepsilon)$ such that:
    (a) $\Phi^*(\cdot, y_0)$ is upper semicontinuous at $y_0$ on $E_0$.
    (b) $\Phi^*(y_0, y_0) \leq \Phi^*(x_0, x_0) + \varepsilon$.
    Then for every $A \in \mathcal{E}$,
    $$\liminf_n a_n^{-1} \log \mathbf{P}\{Y_n \in A\} \geq - \inf_{x \in A^\circ} \Phi^*(x, x).$$
    Moreover, the level sets $\{x \in E : \Phi^*(x, x) \leq l\}, l \geq 0$, are compact.

REMARK 2.3. It is easily shown that, in the presence of the first part of assumption 9, the condition "$x_n \to x$ in $E_0$ implies $\Phi(x_n, \cdot)$ converges to $\Phi(x, \cdot)$ uniformly over the balls in $F$" suffices for 8.

We will need the following two preliminary results.

LEMMA 2.4. *Let $V$ be a separable Banach space. Let $\{\phi_j\}_{j \in N}, \phi$ be proper $w^*$-lower semicontinuous convex functions on $V^*$, and assume:*

(i) $\phi_j(0) = \phi(0) = 0$ *for all $j$.*
(ii) *If $\xi_j \xrightarrow[w^*]{} \xi$, then $\phi(\xi) \leq \liminf_j \phi_j(\xi_j)$.*

*Then, for every $x \in V$, there exists a sequence $\{x_j\}_{j \in N} \subset V, x_j \to x$ such that $\limsup_j \phi_j^*(x_j) \leq \phi^*(x)$.*

This is a particular case of Theorem 1.2 of [17] (the spaces $E$ and $F$ of [17] are here $E = V^*, F = V$).

Let us recall the definition of dominating point [6]. We state it here in our present Banach space framework.

DEFINITION 2.5. Let $V$ be a Banach space and let $\phi : V^* \to \bar{\mathbf{R}}$ be a convex function. Let $D$ be an open convex subset of $V$ such that $D \cap \operatorname{dom} \phi^*$ is nonempty. A point $x_0 \in E$ is a dominating point for $(D, \phi)$ if:

1. $x_0 \in \partial D$.
2. $\phi^*(x_0) = \inf_{x \in D} \phi^*(x)$.

3. There exists $\xi_0 \in E^*$ such that $D \subset \{x : \langle x, \xi_0 \rangle > \langle x_0, \xi_0 \rangle\}$ and $\phi^*(x_0) = \langle x_0, \xi_0 \rangle - \phi(\xi_0)$.

LEMMA 2.6. *Let $V, \phi$ be as in Definition 2.5. Assume:*

(i) $\phi(0) = 0$.
(ii) *For all $a \geq 0$, $L(\phi^*, a)$ is compact.*
(iii) *$\phi$ is $V$-Gâteaux differentiable on $V^*$ and, for all $\xi \in V^*$, $\phi'_\xi$ is continuous, where for $t \in \mathbf{R}$,*

$$\phi_\xi(t) \triangleq \phi(t\xi).$$

*Let $D$ be an open convex subset of $V$ such that $D \cap \operatorname{dom} \phi^*$ is nonempty and $\inf_{x \in D} \phi^*(x) > 0$. Then:*

(a) *There is a unique point $x_0$ satisfying 1 and 2 of Definition 2.5.*
(b) *There exists $\xi_0 \in V^*$ such that $\xi_0$ satisfies 3 of Definition 2.5 and*

$$\nabla \phi(\xi_0) = x_0.$$

This is a particular case of Theorem 2.3 of [6]. A point that should be emphasized is that under the assumptions of Lemma 2.6, the unique point $x_0$ in (a) is automatically a dominating point.

PROOF OF THEOREM 2.2(a). The key part of the proof is to show that under the present hypotheses, condition (11) of Theorem 2.2 of [5] holds. Let $x_0 \in E$ be such that $\Phi^*(x_0, x_0) < \infty$. We must show: for every $\varepsilon > 0$, there exists $x_1 \in E$ such that $x_1 \in B(x_0, \varepsilon), \partial \Phi^*(x_1, x_1) \neq \phi$ and

(2.1) $$\Phi^*(x_1, x_1) < \Phi^*(x_0, x_0) + \varepsilon.$$

Let $y_0$ be as in assumption 11 of the present theorem. Then by assumptions 7 and 11, there exists $\delta_0 > 0$ such that

(2.2) $$\Phi^*(x, y_0) \leq \Phi^*(y_0, y_0) + 1 \quad \text{for } x \in B(y_0, \delta_0),$$

(2.3) $$\bigcup \{L(\Phi^*(y, \cdot), a) : y \in B(y_0, \delta_0)\}$$

is a relatively compact subset of $E_0$ for all $a \geq 0$.

For $0 < \delta < \delta_0, x \in B(y_0, \delta_0)$, let

$$\alpha(x, \delta) = \inf\{\Phi^*(x, y) : y \in B(y_0, \delta)\},$$
$$\beta(\delta) = \sup\{\alpha(x, \delta) : x \in B(y_0, \delta)\}.$$

It follows from (2.2) that $\beta(\delta) \leq \Phi^*(y_0, y_0) + 1$. For $0 < \delta < \delta_0$, let

$$K_\delta = \operatorname{co}\left(\bigcup \{L(\Phi^*(x, \cdot), \beta(\delta)) \cap \overline{B(y_0, \delta)} : x \in \overline{B(y_0, \delta)}\}\right),$$

where co $A$ is the closed convex hull of $A \subset E$. Then by (2.3) and the fact that $E_0$ is a Banach space, $K_\delta$ is a compact convex subset of $E_0$. We define the map $\rho_\delta : K_\delta \to K_\delta$ by

$$\rho_\delta(x) = \begin{cases} \text{dominating point for } (B(y_0,\delta), \Phi(x,\cdot)), & \text{if } \alpha(x,\delta) > 0, \\ \nabla \Phi(x,0), & \text{if } \alpha(x,\delta) = 0. \end{cases}$$

Then for all $x \in K_\delta, z \in B(y_0, \delta)$,

(2.4) $\quad\quad\quad\quad \rho_\delta(x) \in K_\delta \quad \text{and} \quad \Phi^*(x, \rho_\delta(x)) \leq \Phi^*(x, z).$

In the first case $\rho_\delta(x)$ exists and is unique on account of Lemma 2.6, assumption 9 and (2.3), and (2.4) follows from Definition 2.5 and the definition of $K_\delta$.

In the second case, note first that $\Phi^*(x, \nabla \Phi(x,0)) = 0$. Also, by the compactness of $L(\Phi^*(x,\cdot), 1)$, there exists $y \in \overline{B(y_0, \delta)}$ such that $\Phi^*(x,y) = 0$. By assumption 9 and Lemma 2.4 of [6], we must have $y = \nabla \Phi(x,0)$, and therefore $\rho_\delta(x)$ satisfies (2.4).

We claim now:

(2.5) $\quad\quad\quad\quad\quad\quad\quad\quad \rho_\delta \text{ is continuous.}$

Let $x(n)(n \in \mathbf{N}), x \in K_\delta, x(n) \to x$. Given a subsequence $\{n_k, k \in \mathbf{N}\}$ of $\mathbf{N}$, by the compactness of $K_\delta$ there is a subsequence $\{x(n_{k_j}), j \in \mathbf{N}\}$ of $\{x(n_k), k \in \mathbf{N}\}$ and a point $y \in K_\delta$ such that

$$\rho_\delta(x(n_{k_j})) \to y.$$

Since $\Phi^*$ is jointly semicontinuous, we have

(2.6) $\quad\quad\quad\quad \Phi^*(x,y) \leq \liminf_j \Phi^*(x(n_{k_j}), \rho_\delta(x(n_{k_j}))).$

Let $z \in B(y_0, \delta)$. By Lemma 2.4 with $V = E_0, V^* = F, \phi_j = \Phi(x(n_{k_j}), \cdot), \phi = \Phi(x, \cdot)$ and by assumptions 8 and 9 (note that the latter implies the $w^*$-lower semicontinuity of $\phi_j$ and $\phi$), there exists a sequence $\{z_j, j \in \mathbf{N}\}$ which converges to $z$ and satisfies

(2.7) $\quad\quad\quad\quad \limsup_j \Phi^*(x(n_{k_j}), z_j) \leq \Phi^*(x,z).$

Since $z_j \in B(y_0, \delta)$ for sufficiently large $j$, we have by (2.4), (2.6) and (2.7),

$$\Phi^*(x,y) \leq \Phi^*(x,z) \quad \text{for all } z \in B(y_0, \delta).$$

If $\Phi^*(x,y) > 0$, then $\Phi^*(x(n_{k_j}), \rho_\delta(x(n_{k_j}))) > 0$ for sufficiently large $j$ and $\rho_\delta(x(n_{k_j})) \in \partial B(y_0, \delta)$ (see Remark 2.2(2) of [6]), and it follows that $y \in \partial B(y_0, \delta)$. Therefore $y$ is the dominating point for $(B(y_0, \delta), \Phi(x, \cdot))$; that is, $y = \rho_\delta(x)$. If $\Phi^*(x,y) = 0$, then $y = \nabla \Phi(x,0) = \rho_\delta(x)$ by Lemma 2.4 of [6].

We have shown: for every subsequence $\{n_k, k \in \mathbf{N}\}$ of $\mathbf{N}$, there exists a subsequence $\{\rho_\delta(x(n_{k_j})), j \in \mathbf{N}\}$ of $\{\rho_\delta(x(n_k)), k \in \mathbf{N}\}$ such that $\rho_\delta(x(n_{k_j})) \to \rho_\delta(x)$. This proves (2.5).

By the Schauder–Tychonoff fixed point theorem (see, e.g., [14], page 143), there exists $y_\delta \in K_\delta$ such that $\rho_\delta(y_\delta) = y_\delta$. By Lemma 2.6, there exists $\xi_\delta \in F$ such that
$$\nabla \Phi(y_\delta, \xi_\delta) = y_\delta,$$
which implies $\xi_\delta \in \partial \Phi^*(y_\delta, y_\delta)$. Since $y_\delta = \rho_\delta(y_\delta) \in \overline{B(y_0, \delta)}$, we have $y_\delta \to y_0$ as $\delta \to 0$, and by assumption 11,
$$\limsup_\delta \Phi^*(y_\delta, y_\delta) \leq \limsup_\delta \Phi^*(y_\delta, y_0) \leq \Phi^*(y_0, y_0).$$
Taking now $x_1 = y_\delta$ for sufficiently small $\delta$, (2.1) is satisfied.

(b) We will now show that under (2.1) and assumptions 1–6, 9 and 10, the proof of Theorem 2.2 of [5] goes through. Proceeding as in [5], page 490, we have: for $K$ compact,
$$\sup_{y \in V \cap K} [\langle y, a_n \xi \rangle - \Phi_n(y, a_n \xi)] \leq a_n(b_n(K, \xi) + \Phi^*(x_0, x_0) + \varepsilon).$$
Therefore,
$$\mathbf{P}_n\{Y_n \in A\} \geq \mathbf{P}_n\{Y_n \in V \cap K\}$$
$$\geq \exp[-a_n(b_n(K, \xi) + \Phi^*(x_0, x_0) + \varepsilon)]$$
$$\times \int I_{V \cap K}(Y_n) \exp[\langle Y_n, a_n \xi \rangle - \Phi_n(Y_n, a_n \xi)] d\mathbf{P}_n.$$
As in [5], in order to obtain the lower bound, it is enough to show that $\lim_n \sup \mathbf{P}_{n,\xi}\{Y_n \in (V \cap K)^c\} = 0$, or

(2.8) $$\limsup_n \mathbf{P}_{n,\xi}\{Y_n \in V^c \cap K\} = 0,$$

(2.9) $$\limsup_n \mathbf{P}_{n,\xi}\{Y_n \in K^c\} = 0.$$

For $y \in E, \eta \in F$, let
$$\Phi_{n,\xi}(y, \eta) = \Phi_n(y, a_n \xi + \eta) - \Phi_n(y, a_n \xi),$$
$$\Phi_\xi(y, \eta) = \Phi(x, \xi + \eta) - \Phi(x, \xi).$$
Then by assumption 4, for all $\eta \in F$,

(2.10) $$\mathbf{E}_{n,\xi} \exp[\langle Y_n, \eta \rangle - \Phi_{n,\xi}(Y_n, \eta)] = 1.$$

From assumption 3, it easily follows that for every compact set $K \subset E, \eta \in F$,
$$\limsup_n \sup_{x \in K} |a_n^{-1} \Phi_{n,\xi}(x, a_n \eta) - \Phi_\xi(x, \eta)| = 0.$$

By Theorem 2.1, for any compact set $K \subset E$,

(2.11) $\quad \limsup_n a_n^{-1} \log \mathbf{P}_{n,\xi}\{Y_n \in V^c \cap K\} \leq -\inf\{\Phi_\xi^*(x,x) : x \in V^c \cap K\}.$

As in [5], condition 10 implies that the expression in the right-hand side of (2.11) is strictly negative, which proves (2.8). Therefore the proof will be complete if we show that (2.9) holds for a suitable choice of $K$. In fact, $\{\mathcal{L}_{P_{n,\xi}}(Y_n)\}_{n \in N}$ is exponentially tight, to wit

$$\begin{aligned}
& \mathbf{P}_{n,\xi}\{Y_n \in K^c\} \\
& = \int I_{K^c}(Y_n) \exp[\langle Y_n, a_n\xi\rangle - \Phi_n(Y_n, a_n\xi)] \, d\mathbf{P}_n \\
(2.12) \quad & \leq (\mathbf{P}_n\{Y_n \in K^c\})^{1/2} \left(\int \exp[\langle Y_n, 2a_n\xi\rangle - 2\Phi_n(Y_n, a_n\xi)] \, d\mathbf{P}_n\right)^{1/2} \\
& \leq (\mathbf{P}_n\{Y_n \in K^c\})^{1/2} \exp\left[\sup_{x \in E} |\Phi_n(x, 2a_n\xi)| + 2\sup_{x \in E} |\Phi_n(x, a_n\xi)|\right],
\end{aligned}$$

and it follows from (2.12) and assumptions 5 and 6 that, given $b > 0$, $K$ may be chosen so that

$$\limsup_n a_n^{-1} \log \mathbf{P}_{n,\xi}\{Y_n \in K^c\} < -b.$$

The compactness of the level sets of the rate function follows from exponential tightness and the lower bound by a well-known argument. □

**3. Application to a stochastic recursive scheme.** Let $\mu: \mathbf{R}^d \times \mathcal{B}(\mathbf{R}^d) \to [0,1]$ be a Markov kernel. We will consider an i.i.d. sequence of random vector fields with Markov kernel $\mu$; that is, let $F_j: \mathbf{R}^d \times \Omega \to \mathbf{R}^d, j \in \mathbf{N}$, be a sequence of measurable maps such that:

(3.1) For all $j \in \mathbf{N}, x \in \mathbf{R}^d, \mathcal{L}(F_j(x)) = \mu(x, \cdot)$ [here $F_j(x) \equiv F_j(x, \cdot)$].
(3.2) If $j, k \in \mathbf{N}, j \neq k$, then $\{F_j(x) : x \in \mathbf{R}^d\}$ and $\{F_k(x) : x \in \mathbf{R}^d\}$ are independent and have the same distribution.

For $\alpha \in \mathbf{R}^d, x \in \mathbf{R}^d$, let $G: \mathbf{R}^d \times \mathbf{R}^d \to \mathbf{R}$ be defined by

$$G(x, \alpha) = \log \int \exp(\langle y, \alpha \rangle) \mu(x, dy).$$

We shall consider the following conditions:

(3.3) For each $\alpha \in \mathbf{R}^d$, $\sup_{x \in \mathbf{R}^d} G(x, \alpha) < \infty$.
(3.4) $G$ is continuous.
(3.5) For each $r > 0$, the family of functions $\{\nabla G(\cdot, \alpha) : \alpha \in \bar{B}(0, r)\}$ satisfies a uniform Lipschitz condition on $\bar{B}(0, r)$. That is, there exists a constant $D(r) > 0$ such that, for all $\alpha, x, y \in \bar{B}(0, r)$,

$$\|\nabla G(y, \alpha) - \nabla G(x, \alpha)\| \leq D(r)\|y - x\|,$$

where $\|\cdot\|$ is the Euclidean norm on $\mathbf{R}^d$, $\bar{B}(0,r) = \{x \in \mathbf{R}^d : \|x\| \leq r\}$ and $\triangledown G(y,\alpha)$ is the gradient of $G(y,\cdot)$, evaluated at $\alpha$ [$\triangledown G(y,\alpha)$ exists by (3.3)].

(3.6) For every $\alpha \in \mathbf{R}^d$, every $b > 0$, there exists $\tau > 0$ such that

$$\sup\{\mathbf{E}\exp[\tau(\|y-z\|)^{-1}\langle F_1(y) - F_1(z), \alpha\rangle] : y \neq z, \|y\| \leq b, \|z\| \leq b\} < \infty.$$

For fixed $x \in \mathbf{R}^d, n \in \mathbf{N}$, we define recursively, for $0 \leq k \leq n$, the $\mathbf{R}^d$-valued r.v.'s

(3.7)
$$X^x_{n,0} = x,$$
$$X^x_{n,k} = X^x_{n,k-1} + n^{-1} F_k(X^x_{n,k-1}), \qquad k \geq 1,$$

so that

(3.8) $$X^x_{n,k} = x + n^{-1} \sum_{j=1}^{k} F_j(X^x_{n,j-1}), \qquad k \geq 1.$$

Let $T = [0,1]$, and let $\{Y^x_n\}_{n \in N}$ be the $C(T, \mathbf{R}^d)$-valued random vectors given by

$$Y^x_n(t) = \begin{cases} X^x_{n,k}, & \text{if } t = k/n, k = 0, \ldots, n, \\ \text{defined by linear interpolation on } \left[\dfrac{k-1}{n}, \dfrac{k}{n}\right], & k = 1, \ldots, n \end{cases}$$

$$= X^x_{n,[nt]} + (nt - [nt])(n^{-1} F_{[nt]+1}(X^x_{n,[nt]})),$$

where $[\cdot]$ is the integer part function.

For $y, z \in \mathbf{R}^d$, let

$$G^*(y,z) = \sup_{\alpha \in R^d}[\langle z, \alpha\rangle - G(y,\alpha)].$$

Let $\mathcal{C}$ be the Borel $\sigma$-algebra of $C(T, \mathbf{R}^d)$.

THEOREM 3.1. *Assume* (3.3)–(3.6). *Then* $\{\mathcal{L}(Y^x_n)\}_{n \in N}$ *satisfies the large deviation principle on* $C(T, \mathbf{R}^d)$ *with rate function*

$$I^x(f) = \begin{cases} \displaystyle\int_T G^*(f(s), f'(s))\, ds, & \text{if } f(0) = x \\ & \text{and } f \text{ is absolutely continuous,} \\ \infty, & \text{otherwise.} \end{cases}$$

*More specifically, under conditions* (3.3) *and* (3.4), *the upper bound holds:*

for all $A \in \mathcal{C}$, $\quad \limsup_n n^{-1} \log \mathbf{P}\{Y^x_n \in A\} \leq -\inf\{I^x(f) : f \in \bar{A}\},$

*and under conditions* (3.3)–(3.6) *the lower bound holds:*

for all $A \in \mathcal{C}$, $\quad \liminf_n n^{-1} \log \mathbf{P}\{Y^x_n \in A\} \geq -\inf\{I^x(f) : f \in A^0\}.$

*Moreover, under conditions* (3.3) *and* (3.4), *the level sets*

$$L(I^x, \ell) \triangleq \{f \in C(T, \mathbf{R}^d) : I^x(f) \leq \ell\}(\ell \geq 0) \text{ are compact.}$$

REMARK 3.2. 1. The large deviation principle for $\{\mathcal{L}(Y_n^x)\}_{n\in N}$ is presented in [7], Theorem 6.3.3, by a different approach, under condition (3.3), an assumption that together with (3.3) implies (3.4), and either (i) an assumption on the supports of the measures $\mu(x,\cdot)$ or (ii) a special Lipschitz-type assumption on $G^*$. The relation between (i) or (ii) and conditions (3.5), (3.6) is not immediately clear.

2. Conditions (3.3)–(3.6) are hypotheses on the data of the problem—in the sense that (3.3)–(3.5) are assumptions on the Laplace transforms of $\{F(y):y\in\mathbf{R}^d\}$ and (3.6) is an assumption on the Laplace transforms of $\{F(y)-F(z):y,z\in\mathbf{R}^d\}$—and not on $G^*$. Condition (3.5) is used only to verify condition 10 of Theorem 2.2. Condition (3.6) is used only in Lemma 3.7.

3. In the broad class of cases presented later in Theorem 3.11, conditions (3.3)–(3.6) are quite easy to verify.

For the proof of Theorem 3.1, we need several lemmas. The first one, which is elementary, gives a useful expression for functions on $T$ defined by linear interpolation.

LEMMA 3.3. *For $n\in\mathbf{N}, i=1,\ldots,n$, let*
$$\varphi_{ni}(t) = (nt-(i-1))I_{[(i-1)/n,i/n]}(t) + I_{[i/n,1]}(t), \qquad t\in T.$$
*Given $a_i\in\mathbf{R}^d, i=0,\ldots,n$, let*
$$f(t) = \begin{cases} a_i, & \text{for } t=i/n, i=0,\ldots,n, \\ \text{defined by linear interpolation on } \left[\dfrac{i-1}{n}, \dfrac{i}{n}\right], & i=1,\ldots,n. \end{cases}$$
*Then for all $t\in T$,*

$$f(t) = a_0 + \sum_{i=1}^n (a_i - a_{i-1})\varphi_{ni}(t). \tag{3.9}$$

PROOF. For $t\in[\frac{j-1}{n},\frac{j}{n}), j=1,\ldots,n$, we have
$$\varphi_{ni}(t) = \begin{cases} 1, & \text{if } i<j, \\ (nt-(j-1)), & \text{if } i=j, \\ 0, & \text{if } i>j. \end{cases}$$
Therefore, if $g(t)$ is the right-hand side of (3.9), we have for $t\in[\frac{j-1}{n},\frac{j}{n}), j=1,\ldots,n$,
$$g(t) = a_0 + \sum_{i=1}^{j-1}(a_i - a_{i-1}) + (a_j - a_{j-1})(nt-(j-1))$$
$$= a_{j-1} + (a_j - a_{j-1})(nt-(j-1))$$

and $g(1) = a_n$. But this is precisely the definition of $f$. □

It will be convenient for the proof of the lower bound to introduce a perturbation of $\{Y_n^x\}$, as follows. Let $a \geq 0$ and let $\{G_j\}_{j \in N}$ be an independent sequence of $\mathbf{R}^d$-valued r.v.'s with $\mathcal{L}(G_j) = \gamma_d (j \in \mathbf{N})$, where $\gamma_d$ is the canonical Gaussian measure on $\mathbf{R}^d$. We assume also that $\{G_j\}_{j \in N}$ and $\{F_j(x) : j \in \mathbf{N}, x \in \mathbf{R}^d\}$ are independent.

For fixed $x \in \mathbf{R}^d, n \in \mathbf{N}$, define recursively for $0 \leq k \leq n$ the $\mathbf{R}^d$-valued r.v.'s

$$X_{n,0}^{x,a} = x,$$
$$X_{n,k}^{x,a} = X_{n,k-1}^{x,a} + n^{-1}(F_k(X_{n,k-1}^{x,a}) + aG_k), \qquad k \geq 1,$$

so that

$$X_{n,k}^{x,a} = x + n^{-1}\left(\sum_{j=1}^k [F_j(X_{n,j-1}^{x,a}) + aG_j]\right), \qquad k \geq 1.$$

Let $Y_n^{x,a}, n \in \mathbf{N}$, be the $C(T, \mathbf{R}^d)$-valued random vectors given by

$$Y_n^{x,a}(t) = \begin{cases} X_{n,k}^{x,a}, & \text{if } t = k/n, k = 0, \ldots, n, \\ \text{defined by linear interpolation on } \left[\dfrac{k-1}{n}, \dfrac{k}{n}\right], & k = 1, \ldots, n, \end{cases}$$
$$= X_{n,[nt]}^{x,a} + (nt - [nt])(n^{-1}[F_{[nt]+1}(X_{n,[nt]}^{x,a}) + aG_{[nt]+1}]).$$

Let $M(T, \mathbf{R}^d)$ be the space of $\mathbf{R}^d$-valued vector measures defined on the Borel $\sigma$-algebra of $T$. For $f \in C(T, \mathbf{R}^d), \lambda \in M(T, \mathbf{R}^d)$, let

$$\langle f, \lambda \rangle \triangleq \int_T \langle f, d\lambda \rangle.$$

LEMMA 3.4. *For $x \in \mathbf{R}^d, f \in C(T, \mathbf{R}^d), \lambda \in M(T, \mathbf{R}^d)$, let*

$$\Phi_n^{x,a}(f, \lambda) = \langle x, \lambda(T) \rangle + \sum_{i=1}^n G^a\left(f\left(\frac{i-1}{n}\right), n^{-1}\int \varphi_{ni} \, d\lambda\right),$$

*where $\{\varphi_{ni}\}$ are as in Lemma 3.3 and for $y, \alpha \in \mathbf{R}^d$,*

(3.10) $$G^a(y, \alpha) = G(y, \alpha) + \frac{a^2}{2}\|\alpha\|^2.$$

*Then for all $n \in \mathbf{N}$, $\lambda \in M(T, \mathbf{R}^d), a \geq 0$,*

$$\mathbf{E}\exp[\langle Y_n^{x,a}, \lambda\rangle - \Phi_n^{x,a}(Y_n^{x,a}, \lambda)] = 1.$$

PROOF. By Lemma 3.3, we can write, for $t \in T$,

$$Y_n^{x,a}(t) = x + \sum_{i=1}^{n} Z_{ni} \varphi_{ni}(t), \tag{3.11}$$

where $Z_{ni} = n^{-1}[F_i(X_{n,i-1}^{x,a}) + aG_i]$.

Therefore

$$\langle Y_n^{x,a}, \lambda \rangle = \langle x, \lambda(T) \rangle + \sum_{i=1}^{n} \left\langle Z_{ni}, \int \varphi_{ni} \, d\lambda \right\rangle.$$

Let

$$U_k = \sum_{i=1}^{k} \left\langle Z_{ni}, \int \varphi_{ni} \, d\lambda \right\rangle,$$

$$V_k = \sum_{i=1}^{k} G^a \left( X_{n,i-1}^{x,a}, n^{-1} \int \varphi_{ni} \, d\lambda \right);$$

then

$$\langle Y_n^{x,a}, \lambda \rangle - \Phi_n^{x,a}(Y_n^{x,a}, \lambda) = U_n - V_n.$$

Let $\mathcal{F}_k = \sigma(\{F_j(y) : j \leq k, y \in \mathbf{R}^d\} \cup \{G_j : j \leq k\}), k \geq 1$. Then

$$\begin{aligned}
\mathbf{E} \exp[U_n - V_n] \\
&= \mathbf{E}\mathbf{E}(\exp[U_n - V_n] | \mathcal{F}_{n-1}) \\
&= \mathbf{E}\left\{ \exp\left(U_{n-1} - V_{n-1} - G^a\left(X_{n,n-1}^{x,a}, n^{-1} \int \varphi_{nn} \, d\lambda\right)\right) \right. \\
&\qquad \left. \times \mathbf{E}\left[ \exp\left\langle Z_{nn}, \int \varphi_{nn} \, d\lambda \right\rangle \middle| \mathcal{F}_{n-1} \right] \right\},
\end{aligned} \tag{3.12}$$

since $U_{n-1}, V_{n-1}$ and $X_{n,n-1}^{x,a}$ are $\mathcal{F}_{n-1}$-measurable. Next, since $\{F_n(y) : y \in \mathbf{R}^d\} \cup \{G_n\}$ is independent of $\mathcal{F}_{n-1}$ and $\{F_n(y) : y \in \mathbf{R}^d\}$ is independent of $G_n$, we have

$$\begin{aligned}
\mathbf{E}\left[ \exp\left\langle Z_{nn}, \int \varphi_{nn} \, d\lambda \right\rangle \middle| \mathcal{F}_{n-1} \right] \\
&= g(X_{n,n-1}^{x,a}) \\
&= \exp\left[ G\left( X_{n,n-1}^{x,a}, n^{-1} \int \varphi_{nn} \, d\lambda \right) + \frac{a^2 n^{-2}}{2} \left\| \int \varphi_{nn} \, d\lambda \right\|^2 \right] \\
&= \exp\left[ G^a\left( X_{n,n-1}^{x,a}, n^{-1} \int \varphi_{nn} \, d\lambda \right) \right],
\end{aligned} \tag{3.13}$$

where

$$g(y) = \mathbf{E} \exp\left\langle F_n(y) + aG_n, n^{-1} \int \varphi_{nn} \, d\lambda \right\rangle.$$

By (3.12) and (3.13),
$$\mathbf{E}\exp[U_n - V_n] = \mathbf{E}\exp[U_{n-1} - V_{n-1}].$$

Iterating this procedure, we obtain

$$\mathbf{E}\exp[U_n - V_n]$$
$$= \mathbf{E}\exp[U_1 - V_1]$$
$$= \exp\left[-G^a\left(x, n^{-1}\int \varphi_{n1}\,d\lambda\right)\right]\mathbf{E}\exp\left\langle F_1(x) + aG_1, n^{-1}\int \varphi_{n1}\,d\lambda\right\rangle$$
$$= 1. \qquad \square$$

LEMMA 3.5. *For $x \in \mathbf{R}^d, f \in C(T, \mathbf{R}^d), \lambda \in M(T, \mathbf{R}^d)$, let*
$$\Phi^{x,a}(f, \lambda) = \langle x, \lambda(T)\rangle + \int_T G^a(f(s), \lambda([s,1]))\,ds,$$

*where $G^a$ is given by (3.10). Then for every compact set $K \subset C(T, \mathbf{R}^d)$, every $\lambda \in M(T, \mathbf{R}^d)$,*
$$\lim_n \sup_{f \in K} |n^{-1}\Phi_n^{x,a}(f, n\lambda) - \Phi^{x,a}(f, \lambda)| = 0.$$

PROOF. For all $f \in C(T, \mathbf{R}^d), \lambda \in M(T, \mathbf{R}^d)$,

$$|n^{-1}\Phi_n^{x,a}(f, n\lambda) - \Phi^{x,a}(f, \lambda)|$$
(3.14)
$$= \left|n^{-1}\sum_{j=1}^n G^a\left(f\left(\frac{j-1}{n}\right), \int \varphi_{nj}\,d\lambda\right) - \int_T G^a(f(s), \lambda([s,1]))\,ds\right|$$
$$\leq \int_T |G^a(f_n(s), \varphi_n(s)) - G^a(f(s), \varphi(s))|\,ds,$$

where
$$f_n(s) = \sum_{j=1}^n f\left(\frac{j-1}{n}\right)I_{[(j-1)/n, j/n)}(s) + f(i)\,I_{\{i\}}(s),$$
$$\varphi_n(s) = \sum_{j=1}^n \left(\int \varphi_{nj}\,d\lambda\right)I_{[(j-1)/n, j/n)}(s),$$

and $\varphi(s) = \lambda([s,1])$. Clearly,

(3.15) $$\|f_n - f\|_\infty \leq w(f, n^{-1}),$$

where $w$ is the usual modulus of continuity: for $g \in C(T, \mathbf{R}^d), \delta > 0$,

(3.16) $$w(g, \delta) = \sup\{\|g(t) - g(s)\| : s, t \in T, |t - s| \leq \delta\}.$$

Since $\varphi_n(s) \to \varphi(s)$ except possibly at countably many points of $T$, by Egoroff's theorem there exists a measurable set $A \subset T$ such that $m(A) < (4c_2)^{-1}\varepsilon$ and $\varphi_n$ converges to $\varphi$ uniformly on $A^c$, where

$$c_2 = \sup\{|G^a(y,\alpha)| : \|y\| \leq c_1, \|\alpha\| \leq \|\lambda\|_v\},$$
$$c_1 = \sup\{\|f\|_\infty : f \in K\},$$

and $\|\cdot\|_v$ is the total variation norm on $M(T, \mathbf{R}^d)$.

By condition (3.4), $G$ is uniformly continuous on $\bar{B}(0, c_1) \times \bar{B}(0, \|\lambda\|_v)$. Therefore there exists $\delta > 0$ such that $y, z \in \bar{B}(0, c_1), \alpha, \beta \in \bar{B}(0, \|\lambda\|_v), \|y - z\| \leq \delta, \|\alpha - \beta\| \leq \delta$ imply

$$(3.17) \qquad |G^a(y,\alpha) - G^a(z,\beta)| < \varepsilon/2.$$

Let $n_0 \in \mathbf{N}$ be such that:

(i) $\sup_{f \in K} w(f, n_0^{-1}) < \delta$,
(ii) $\sup_{s \in A^c} \|\varphi_n(s) - \varphi(s)\| < \delta$ for $n \geq n_0$.

Then by (3.14), (3.15), (3.17), (i) and (ii), for $n \geq n_0, f \in K$, we have

$$\sup_{s \in A^c} |G^a(f_n(s), \varphi_n(s)) - G^a(f(s), \varphi(s))| \leq \varepsilon/2,$$

and therefore

$$\sup_{f \in K} |n^{-1}\Phi_n^{x,a}(f, n\lambda) - \Phi^{x,a}(f, \lambda)| \leq 2c_2 m(A) + \varepsilon/2$$
$$= \varepsilon. \qquad \square$$

LEMMA 3.6. $\{\mathcal{L}(Y_n^{x,a})\}_{n \in N}$ *is exponentially tight.*

PROOF. We first observe that it is enough to show: for every $b > 0, \varepsilon > 0$, there exist $\delta > 0, n_0 \in \mathbf{N}$, such that

$$(3.18) \qquad \mathbf{P}\{w(Y_n^{x,a}, \delta) > \varepsilon\} \leq e^{-bn} \qquad \text{for } n \geq n_0.$$

To justify this claim, we start by noting that in (3.18) one can take $n_0 = 1$. Suppose (3.18) holds. Since for any $f \in C(T, \mathbf{R}^d)$ we have

$$\lim_{\rho \to 0} w(f, \rho) = 0,$$

one can choose $\rho > 0$ so that

$$\mathbf{P}\{w(Y_n^{x,a}, \rho) > \varepsilon\} \leq e^{-bn} \qquad \text{for } n < n_0.$$

Replacing $\delta$ by $\min\{\delta, \rho\}$, (3.18) is now valid for $n \geq 1$.

Next, given $b > 0$, choose $b_j > b (j \geq 1)$ such that

$$\sum_{j=1}^\infty \exp[-(b_j - b)] < 1$$

and let $\varepsilon_j \downarrow 0$. For $\delta_j$ associated to $b_j, \varepsilon_j$ as in (3.18), let

$$K = \{f \in C(T, \mathbf{R}^d) : f(0) = x \text{ and for all } j \in \mathbf{N}, w(f, \delta_j) \leq \varepsilon_j\}.$$

By the Arzelá–Ascoli theorem (see, e.g., [2], page 221), $K$ is compact. For all $n \geq 1$,

$$\mathbf{P}\{Y_n^{x,a} \in K^c\} \leq \sum_{j=1}^{\infty} \mathbf{P}\{w(Y_n^{x,a}, \delta_j) > \varepsilon_j\}$$

$$\leq \sum_{j=1}^{\infty} e^{-b_j n}$$

$$\leq e^{-bn},$$

which proves that (3.18) implies that $\{\mathcal{L}(Y_n^{x,a})\}_{n \in N}$ is exponentially tight.

Let $B = \{\alpha_1, \ldots, \alpha_d\}$ be a basis of $\mathbf{R}^d$ such that, for all $v \in \mathbf{R}^d$,

$$\|v\| \leq \sup_{1 \leq j \leq d} |\langle v, \alpha_j \rangle|.$$

Then, for all $v \in \mathbf{R}^d$,

$$\exp(\|v\|) \leq \sum_{j=1}^{d} [\exp(\langle v, \alpha_j \rangle) + \exp(\langle v, -\alpha_j \rangle)].$$

Therefore, for all $y \in \mathbf{R}^d, \tau > 0$,

$$(3.19) \quad \mathbf{E} \exp(\tau \|F_1(y)\|) \leq \sum_{j=1}^{d} [\exp(G(y, \tau \alpha_j)) + \exp(G(y, -\tau \alpha_j))]$$

$$\leq 2d \max\{\exp(G(y, \tau \alpha)) : \alpha \in B \cup (-B)\}.$$

Using condition (3.3), it follows that

$$c_1 = c_1(\tau) = \sup_{y \in R^d} \mathbf{E} \exp(\tau \|F_1(y)\|) < \infty.$$

We claim next that if $c_2 = c_2(\tau) = \mathbf{E} \exp(\tau a \|G_1\|)$, then for $p, q \in \mathbf{N}, 1 \leq p < q \leq n$,

$$(3.20) \quad \mathbf{E} \exp\left(\sum_{j=p+1}^{q} \tau \|F_j(X_{n,j-1}^{x,a}) + aG_j\|\right) \leq (c_1 c_2)^{q-p}.$$

Arguing similarly to the proof of Lemma 3.4,

$$\mathbf{E}\left(\sum_{j=p+1}^{q} \tau \|F_j(X_{n,j-1}^{x,a}) + aG_j\|\right)$$

$$= \mathbf{E}\bigg(\exp\bigg(\sum_{j=p+1}^{q-1} \tau\|F_j(X_{n,j-1}^{x,a}) + aG_j\|\bigg)$$
$$\times \mathbf{E}[\exp(\tau\|F_q(X_{n,q-1}^{x,a}) + aG_q\|)|\mathcal{F}_{q-1}]\bigg).$$

But

$$\mathbf{E}[\exp(\tau\|F_q(X_{n,q-1}^{x,a}) + aG_q\|)|\mathcal{F}_{q-1}]$$
$$\leq \mathbf{E}[\exp(\tau\|F_q(X_{n,q-1}^{x,a})\|)|\mathcal{F}_{q-1}]\mathbf{E}\exp(\tau a\|G_q\|)$$
$$= g(X_{n,q-1}^{x,a}) \cdot c_2(\tau)$$
$$\leq c_1(\tau)c_2(\tau),$$

where $g(y) = \mathbf{E}\exp(\tau\|F_q(y)\|)$. Now (3.20) follows by iteration.

The next step is to show that, for $m \in \mathbf{N}, m \leq n$,

$$(3.21) \qquad w(Y_n^{x,a}, m^{-1}) \leq 3 \sup_{0 \leq i \leq m-1} \sum_{j=[(ni)/m]+1}^{[(n(i+1))/m]+1} \|Z_{nj}\|,$$

where $Z_{nj}$ is as in Lemma 3.4. First we note that by the triangle inequality, for any $f \in C(T, \mathbf{R}^d)$,

$$(3.22) \qquad w(f, m^{-1}) \leq 3 \sup_{0 \leq i \leq m-1} \sup_{t \in [i/m, (i+1)/m]} \left\|f(t) - f\left(\frac{i}{m}\right)\right\|.$$

For $t \in [\frac{i}{m}, \frac{i+1}{m}]$,

$$(3.23) \qquad \begin{aligned} &\left\|Y_n^{x,a}(t) - Y_n^{x,a}\left(\frac{i}{m}\right)\right\| \\ &= \bigg\|\bigg(x + \sum_{j=1}^{[nt]} Z_{nj} + (nt - [nt])Z_{n,[nt]+1}\bigg) \\ &\quad - \bigg(x + \sum_{j=1}^{[(ni)/m]} Z_{nj} + \bigg(\frac{ni}{m} - \bigg[\frac{ni}{m}\bigg]\bigg)Z_{n,[(ni)/m]+1}\bigg)\bigg\| \\ &\leq \sum_{j=[(ni)/m]+1}^{[(n(i+1))/m]+1} \|Z_{nj}\|. \end{aligned}$$

Now (3.21) follows from (3.22) and (3.23). For $\varepsilon > 0, \tau > 0$, by (3.20) and (3.21),

$$\mathbf{P}\{w(Y_n^{x,a}, m^{-1}) > \varepsilon\}$$
$$\leq \sum_{i=0}^{m-1} \mathbf{P}\bigg\{\sum_{j=[(ni)/m]+1}^{[(n(i+1))/m]+1} \|Z_{nj}\| > \frac{\varepsilon}{3}\bigg\}$$

$$\leq \sum_{i=0}^{m-1} e^{-\tau n\varepsilon/3} \mathbf{E} \exp\left[\tau \sum_{j=[(ni)/m]+1}^{[(n(i+1))/m]+1} \|F_j(X_{n,j-1}^{x,a}) + aG_j\|\right]$$

$$\leq m e^{-\tau n\varepsilon/3} (c_1(\tau)c_2(\tau))^{(n/m)+2}$$

$$< \exp\left[-n\left(\tau\varepsilon/3 - \frac{m}{n} - \left(\frac{1}{m} + \frac{2}{n}\right)\log(c_1(\tau)c_2(\tau))\right)\right].$$

Given $b > 1$, choose $\tau > 6b\varepsilon^{-1}, m \geq \log(c_1(\tau)c_2(\tau))$. Then

$$\limsup_n n^{-1} \log \mathbf{P}\{w(Y_n^{x,a}, m^{-1}) > \varepsilon\} \leq -\tau\varepsilon/3 + 1$$

$$< -b,$$

which establishes (3.18). $\square$

In the next lemma we show that $\{Y_n^x\}_{n \in N}$ and $\{Y_n^{x,a}\}_{n \in N}$ are superexponentially close in probability as $a \to 0$.

LEMMA 3.7. *For every $\varepsilon > 0$,*

$$\lim_{a \downarrow 0} \limsup_n n^{-1} \log \mathbf{P}\{\|Y_n^{x,a} - Y_n^x\|_\infty > \varepsilon\} = -\infty.$$

PROOF. We will use the following estimate: for all $n \in \mathbf{N}$,

(3.24)
$$\sup_{k \leq n} \|X_{n,k}^{x,a} - X_{n,k}^x\|$$
$$\leq n^{-1} a \sum_{j=1}^n \left[\|G_j\| \prod_{i=j+1}^n (1 + n^{-1} H_i(X_{n,i-1}^{x,a}, X_{n,i-1}^x))\right],$$

where, for $y, z \in \mathbf{R}^d$,

$$H_j(y, z) = \begin{cases} (\|y - z\|)^{-1} \|F_j(y) - F_j(z)\|, & \text{if } y \neq z, \\ 0, & \text{if } y = z. \end{cases}$$

To prove (3.24), we use the following elementary inequality, which is obtained at once by induction: If $\{a_k\}_{k \in N}, \{b_k\}_{k \in N}, \{c_k\}_{k \in N}$ are nonnegative real numbers such that

$$a_1 \leq c_1, a_k \leq a_{k-1} b_k + c_k \qquad \text{for } k \geq 2,$$

then for all $k \geq 2$,

(3.25)
$$a_k \leq \sum_{j=1}^{k-1} c_j \left(\prod_{i=j+1}^k b_i\right) + c_k.$$

We have, for $1 \leq k \leq n$,

$$X_{n,k}^{x,a} - X_{n,k}^x$$
$$= (X_{n,k-1}^{x,a} - X_{n,k-1}^x) + n^{-1}(F_k(X_{n,k-1}^{x,a}) - F_k(X_{n,k-1}^x)) + n^{-1}aG_k,$$

and therefore,

(3.26)
$$\|X_{n,k}^{x,a} - X_{n,k}^x\|$$
$$\leq \|X_{n,k-1}^{x,a} - X_{n,k-1}^x\|(1 + n^{-1}H_k(X_{n,k-1}^{x,a}, X_{n,k-1}^x)) + n^{-1}a\|G_k\|.$$

Also,

$$X_{n,1}^{x,a} - X_{n,1}^x = (x + n^{-1}F_1(x) + n^{-1}aG_1) - (x + n^{-1}F_1(x))$$
$$= n^{-1}aG_1.$$

Setting $a_k = \|X_{n,k}^{x,a} - X_{n,k}^x\|, b_k = 1 + n^{-1}H_k(X_{n,k-1}^{x,a}, X_{n,k-1}^x), c_k = n^{-1}a\|G_k\|$, (3.24) follows from (3.25) and (3.26).

Using the elementary inequality $1 + x \leq e^x$ ($x \in \mathbf{R}$), (3.24) implies

$$\|Y_n^{x,a} - Y_n^x\|_\infty = \sup_{k \leq n} \|X_{n,k}^{x,a} - X_{n,k}^x\|$$
$$\leq n^{-1}a\left(\sum_{j=1}^n \|G_j\|\right) \exp\left(n^{-1} \sum_{i=1}^n H_i(X_{n,i-1}^{x,a}, X_{n,i-1}^x)\right).$$
(3.27)

For $\tau > 0, b > 0$, let $c(\tau, b) = \sup\{\mathbf{E}\exp(\tau H_1(y,z)) : \|y\| \leq b, \|z\| \leq b\}$.
For $\alpha \in \mathbf{R}^d$, let

$$c(\alpha, \tau, b) = \sup\{\mathbf{E}\exp(\tau(\|y - z\|)^{-1}\langle F_1(y) - F_1(z), \alpha\rangle) : \|y\| \leq b, \|z\| \leq b\},$$

and let

$$\bar{c}(\tau, b) = 2d \max\{c(\alpha, \tau, b) : \alpha \in B \cup (-B)\},$$

where $B$ is as in the proof of Lemma 3.6.

By conditioning as in Lemma 3.4 and iterating, we have

(3.28)
$$\mathbf{E}\left\{I\left(\sup_{k \leq n-1} \|X_{n,k}^{x,a}\| \leq b, \sup_{k \leq n-1} \|X_{n,k}^x\| \leq b\right)\right.$$
$$\left. \times \exp\left[\tau \sum_{i=1}^n H_i(X_{n,i-1}^{x,a}, X_{n,i-1}^x)\right]\right\}$$
$$\leq \mathbf{E}\left\{I\left(\sup_{k \leq n-2} \|X_{n,k}^{x,a}\| \leq b, \sup_{k \leq n-2} \|X_{n,k}^x\| \leq b\right)\right.$$
$$\left. \times \exp\left[\tau \sum_{i=1}^{n-1} H_i(X_{n,i-1}^{x,a}, X_{n,i-1}^x)\right]\right\} c(\tau, b)$$
$$\leq (c(\tau, b))^n$$
$$\leq \bar{c}(\tau, b)^n;$$

the last inequality is proved similarly to (3.19).

For $\varepsilon > 0, b > 0$,

$$
\begin{aligned}
(3.29) \quad & \mathbf{P}\{\|Y_n^{x,a} - Y_n^x\|_\infty > \varepsilon\} \\
& \leq \mathbf{P}\{\|Y_n^{x,a}\|_\infty > b\} + \mathbf{P}\{\|Y_n^x\|_\infty > b\} \\
& \quad + \mathbf{P}\left\{\sup_{k \leq n}\|X_{n,k}^{x,a}\| \leq b, \sup_{k \leq n}\|X_{n,k}^x\| \leq b, \|Y_n^{x,a} - Y_n^x\|_\infty > \varepsilon\right\}.
\end{aligned}
$$

By (3.27), Markov's inequality and (3.28),

$$
\begin{aligned}
(3.30) \quad & \mathbf{P}\left\{\sup_{k \leq n}\|X_{n,k}^{x,a}\| \leq b, \sup_{k \leq n}\|X_{n,k}^x\| \leq b, \|Y_n^{x,a} - Y_n^x\|_\infty > \varepsilon\right\} \\
& \leq \mathbf{P}\left\{\sup_{k \leq n}\|X_{n,k}^{x,a}\| \leq b, \sup_{k \leq n}\|X_{n,k}^x\| \leq b, \right. \\
& \qquad \left. n^{-1}\sum_{i=1}^n H_i(X_{n,i-1}^{x,a}, X_{n,i-1}^x) > r\right\} \\
& \quad + \mathbf{P}\left\{n^{-1}a\sum_{j=1}^n \|G_j\| > \varepsilon e^{-r}\right\} \\
& \leq e^{-n\tau r}(\bar{c}(\tau,b))^n + \exp(-na^{-1}\varepsilon e^{-r})(\mathbf{E}\exp(\|G_1\|))^n.
\end{aligned}
$$

Next, using (3.11) and (3.20), we have

$$
\|Y_n^{x,a}\|_\infty \leq \|x\| + \sum_{i=1}^n n^{-1}\|F_i(X_{n,i-1}^{x,a}) + aG_i\|,
$$

$$
(3.31) \quad \sup_{0 \leq a \leq 1} \mathbf{E}\exp(n\|Y_n^{x,a}\|_\infty) \leq e^{n\|x\|}(c_1(1)c_2(1))^n,
$$

$$
\sup_{0 \leq a \leq 1} \mathbf{P}\{\|Y_n^{x,a}\|_\infty > b\} \leq e^{-nb}e^{n\|x\|}(c_1(1)c_2(1))^n.
$$

By (3.29),

$$
\begin{aligned}
& \limsup_n n^{-1}\log \mathbf{P}\{\|Y_n^{x,a} - Y_n^x\|_\infty > \varepsilon\} \\
& \leq \max\left\{\limsup_n n^{-1}\log \sup_{0 \leq a \leq 1}\mathbf{P}\{\|Y_n^{x,a}\|_\infty > b\}, -(\tau r - \log \bar{c}(\tau,b)), \right. \\
& \qquad \left. -(a^{-1}\varepsilon e^{-r} - \log \mathbf{E}\exp(\|G_1\|))\right\}.
\end{aligned}
$$

Given $\ell > 0$, by (3.31) there exists $b > 0$ such that

$$
\limsup_n n^{-1}\log \sup_{0 \leq a \leq 1}\mathbf{P}\{\|Y_n^{x,a}\|_\infty > b\} < -\ell.
$$

By condition (3.6), there exists $\tau > 0$ such that $\bar{c}(\tau,b) < \infty$. Let $r > 0$ be such that $\tau r - \log \bar{c}(\tau,b) > \ell$. Then

$$
\limsup_{a \downarrow 0}\limsup_n n^{-1}\log \mathbf{P}\{\|Y_n^{x,a} - Y_n^x\|_\infty > \varepsilon\} < -\ell.
$$

Since $\ell$ is arbitrary, this completes the proof. $\square$

LEMMA 3.8. *Let $G^a$ be given by (3.10). Then $(G^a)^*$ is continuous on $\mathbf{R}^d \times \mathbf{R}^d$, where*

$$(G^a)^*(y, z) = \sup_{\alpha \in \mathbf{R}^d} [\langle z, \alpha \rangle - G^a(y, \alpha)].$$

PROOF. We include for completeness the following argument, which is a slight variant of one to be found, for example, in [3], pages 958 and 959.

By Jensen's inequality, for all $y \in \mathbf{R}^d, \alpha \in \mathbf{R}^d$,

$$\begin{aligned}(3.32) \qquad G^a(y, \alpha) &\geq \int \langle z, \alpha \rangle \mu(y, dz) + \frac{a^2}{2} \|\alpha\|^2 \\ &\geq q(\alpha),\end{aligned}$$

where $q(\alpha) = -D\|\alpha\| + \frac{a^2}{2}\|\alpha\|^2$, for a suitable constant $D$ which exists by condition (3.3). Therefore, by an elementary calculation (see, e.g., [3], page 955),

$$(3.33) \qquad (G^a)^*(y, z) \leq q^*(z) = (2a^2)^{-1}(\|z\| + D)^2,$$

so $(G^a)^*$ is everywhere finite. Suppose $(y(n), z(n)) \to (y, z)$ in $\mathbf{R}^d \times \mathbf{R}^d$. For any positive sequence $\varepsilon_k \downarrow 0$ and any subsequence $\{n_k\}$ of $\{n\}$, there exists $\{\alpha_k\}$ in $\mathbf{R}^d$ such that

$$0 \leq (G^a)^*(y(n_k), z(n_k)) \leq \langle z(n_k), \alpha_k \rangle - G^a(y(n_k), \alpha_k) + \varepsilon_k$$
$$\leq (\|z(n_k)\| + D)\|\alpha_k\| - \frac{a^2}{2}\|\alpha_k\|^2 + \varepsilon_k,$$

and hence $\{\alpha_k\}$ is bounded. Therefore, there exist a subsequence $\{\alpha_{k_j}\}$ of $\{\alpha_k\}$ and $\beta \in \mathbf{R}^d$ such that $\lim_j \alpha_{k_j} = \beta$. Then

$$\begin{aligned}\limsup_j (G^a)^*&(y(n_{k_j}), z(n_{k_j})) \\ &\leq \limsup_j [\langle z(n_{k_j}), \alpha_{k_j} \rangle - G^a(y(n_{k_j}), \alpha_{k_j})] \\ &= \langle z, \beta \rangle - G^a(y, \beta) \\ &\leq (G^a)^*(y, z).\end{aligned}$$

By the lower semicontinuity of $(G^a)^*$,

$$\liminf_j (G^a)^*(y(n_{k_j}), z(n_{k_j})) \geq (G^a)^*(y, z),$$

and therefore

$$\lim_j (G^a)^*(y(n_{k_j}), z(n_{k_j})) = (G^a)^*(y, z).$$

This proves the continuity of $(G^a)^*$ at $(y, z)$. □

PROOF OF THEOREM 3.1. (I) *Upper bounds.* In the context of Theorem 2.1, let $a_n = n, E = C(T, \mathbf{R}^d), F = M(T, \mathbf{R}^d)$ and let $\mathcal{E} = \mathcal{C}$. Also let $Y_n = Y_n^x, \Phi_n = \Phi_n^x, \Phi = \Phi^x$, where for $f \in C(T, \mathbf{R}^d), \lambda \in M(T, \mathbf{R}^d)$,

$$\Phi_n^x(f, \lambda) = \langle x, \lambda(T) \rangle + \sum_{i=1}^{n} G\left(f\left(\frac{i-1}{n}\right), n^{-1} \int \varphi_{ni}\, d\lambda\right),$$

$$\Phi^x(f, \lambda) = \langle x, \lambda(T) \rangle + \int_T G(f(s), \lambda([s, 1]))\, ds.$$

Assume (3.3) and (3.4). It is immediate that conditions 1 and 2 of Theorem 2.1 hold. Conditions 3–5 of Theorem 2.1 hold, respectively, by Lemmas 3.5, 3.4 and 3.6 with $a = 0$. Applying Theorem 2.1, for all $A \in \mathcal{C}$,

$$\limsup_n n^{-1} \log \mathbf{P}\{Y_n^x \in A\} \leq -\inf_{f \in \bar{A}} (\Phi^x)^*(f, f).$$

But $(\Phi^x)^*(f, f) = I^x(f)$ for all $f \in E$ by the argument in Theorem 6.1 of [5], which applies easily to the present situation. This completes the proof of the upper bound.

(II) *Compactness of the level sets.* Assume (3.3) and (3.4). We will frame the argument so that it is useful in the proof of the lower bound. Let $\bar{G}(\alpha) = \sup_{y \in R^d} G(y, \alpha)$, and for $\lambda \in F, f \in E$, let

$$\psi(\lambda) = \|\lambda(T)\| \|x\| + \int_T \bar{G}(\lambda([s, 1]))\, ds,$$

$$\psi^*(f) = \sup_{\lambda \in F}[\langle f, \lambda \rangle - \psi(\lambda)].$$

Then for all $h \geq 0$:

(i) $L(\psi^*, h)$ is compact.
(ii) $\bigcup_{f \in E} L((\Phi^x)^*(f, \cdot), h) \subset L(\psi^*, h)$.

Since $L(\psi^*, h)$ is closed, by the Arzelá–Ascoli theorem (see, e.g, [2], page 221), to prove (i) it is enough to show:

(a) $\sup\{\|f(0)\| : f \in L(\psi^*, h)\} < \infty$,
(b) $\lim_{\delta \downarrow 0} \sup\{w(f, \delta) : f \in L(\psi^*, h)\} = 0$,

where $w$ is given by (3.16). We prove only (b); the proof of (a) is similar but simpler. If $f \in L(\psi^*, h)$, then for all $\lambda \in F$,

(3.34) $$\int_T \langle f, d\lambda \rangle \leq \psi(\lambda) + h.$$

Let $s, t \in T, s < t, \rho > 0, \alpha \in \mathbf{R}^d$, and let $\lambda = \rho\alpha(\delta_t - \delta_s)$. Then by (3.34) we have

$$\rho\langle f(t) - f(s), \alpha \rangle \leq \psi(\rho\alpha(\delta_t - \delta_s)) + h$$

$$= \bar{G}(\rho\alpha)(t-s) + h,$$
$$\|f(t) - f(s)\| = \sup\{\langle f(t) - f(s), \alpha\rangle : \|\alpha\| \leq 1\}$$
$$\leq \rho^{-1}\sup\{\bar{G}(\rho\alpha) : \|\alpha\| \leq 1\}(t-s) + \rho^{-1}h,$$

and therefore

(3.35) $\quad \sup\{w(f,\delta) : f \in L(\psi^*, h)\} \leq \rho^{-1}\sup\{\bar{G}(\rho\alpha) : \|\alpha\| \leq 1\}\delta + \rho^{-1}h.$

Since $\bar{G}$ is a finite convex function by (3.3), hence continuous, (b) follows from (3.35).

To prove (ii): for all $f \in E, \lambda \in F, g \in E$,
$$\Phi^x(f, \lambda) \leq \psi(\lambda),$$
$$(\Phi^x)^*(f, g) \geq \psi^*(g),$$

and therefore for all $h \geq 0$,

(3.36) $\qquad\qquad L((\Phi^x)^*(f, \cdot), h) \subset L(\psi^*, h).$

Finally note that (3.36) implies: for all $h \geq 0$,
$$\{f \in E : (\Phi^x)^*(f, f) \leq h\} \subset L(\psi^*, h),$$

which proves the compactness of the level sets of $I^x$. (Of course, this property also follows from Theorem 2.2 once the lower bound has been established.)

(III) *Lower bounds.* First we prove the lower bound for $\{\mathcal{L}(Y_n^{x,a})\}_{n\in N}$. We take $E_0 = E, \Phi_n = \Phi_n^{x,a}, \Phi = \Phi^{x,a}$. Conditions 1–5 of Theorem 2.1 are proved as in (I). Let $\lambda \in F$. Then for all $n \in \mathbf{N}, f \in E$,

$$n^{-1}\Phi_n^{x,a}(f, n\lambda)$$
$$= \langle x, \lambda(T)\rangle + n^{-1}\sum_{i=1}^{n} G^a\left(f\left(\frac{i-1}{n}\right), \int \varphi_{ni}\,d\lambda\right)$$
$$\leq \langle x, \lambda(T)\rangle + n^{-1}\sum_{i=1}^{n}\left[\bar{G}\left(\int \varphi_{ni}\,d\lambda\right) + \frac{a^2}{2}\left\|\int \varphi_{ni}\,d\lambda\right\|^2\right]$$
$$\leq |\langle x, \lambda(T)\rangle| + \sup\left\{\left|\bar{G}(\alpha) + \frac{a^2}{2}\|\alpha\|^2\right| : \|\alpha\| \leq \|\lambda\|_v\right\}$$
$$\triangleq C,$$

which is finite by the continuity of $\bar{G}$, and therefore
$$\sup_{n}\sup_{f \in E}|n^{-1}\Phi_n^{x,a}(f, n\lambda)| \leq C < \infty.$$

This establishes condition 6 of Theorem 2.2. Condition 7 of Theorem 2.2 for $\Phi^{x,a}$ is proved as in (II) above. It is readily seen that for all $f \in E, \Phi^{x,a}(f, \cdot)$ is convex. The fact that $\Phi^{x,a}(f,\cdot)$ is $E$-Gâteaux differentiable, with $E$-Gâteaux derivative at $\lambda \in F$ given by $\nabla \Phi^{x,a}(f,\lambda) = f_\lambda$, where

$$f_\lambda(t) = x + \int_0^t \nabla G^a(f(s), \lambda([s,1])) \, ds, \qquad t \in T,$$

is proved as in Lemma 7.4 of [5]. Moreover, if

$$\phi(t) = \Phi^{x,a}(f, t\lambda),$$

then

$$\phi'(t) = \langle \nabla \Phi^{x,a}(f, t\lambda), \lambda \rangle$$
$$= \langle x, \lambda(T) \rangle + \int_T \left\langle \int_0^u \nabla G^a(f(s), t\lambda([s,1])) \, ds, d\lambda(u) \right\rangle,$$

and therefore $\phi'$ is continuous. This shows that condition 9 of Theorem 2.2 holds. Next, since

$$|\Phi^{x,a}(f,\lambda) - \Phi^{x,a}(g,\lambda)| \leq \int_T |G(f(s), \lambda([s,1])) - G(g(s), \lambda([s,1]))| \, ds,$$

it follows from condition (3.4) and Remark 2.3 that condition 8 of Theorem 2.2 holds. Using condition (3.5), the fact that condition 10 of Theorem 2.2 holds for $\Phi^{x,a}$ is proved by showing that

$$f = \nabla \Phi^{x,a}(f,\lambda), \qquad g = \nabla \Phi^{x,a}(g,\lambda)$$

imply $f = g$ as in [5], page 518.

Let $(\Phi^{x,a})^*(f_0, f_0) < \infty, \varepsilon > 0$. By the proof of Lemma 7.6 of [5], which applies to the present situation by (3.32), we have: there exists $g_0 \in E$ such that $g_0$ is absolutely continuous, $g_0(0) = x, g_0' \in L^\infty(T)$ and:

(i) $g_0 \in B(f_0, \varepsilon)$,
(ii) $(\Phi^{x,a})^*(g_0, g_0) \leq (\Phi^{x,a})^*(f_0, f_0) + \varepsilon$.

Suppose $f_n \to g_0$ in $E$. Since by (3.33), for almost every $s \in T$,

$$(G^a)^*(f_n(s), g_0'(s)) \leq (2a^2)^{-1}(\|g_0'(s)\| + D)^2,$$

by Lemma 3.8 and the dominated convergence theorem we have

$$(\Phi^{x,a})^*(f_n, g_0) = \int_T (G^a)^*(f_n(s), g_0'(s)) \, ds$$
$$\to \int_T (G^a)^*(g_0(s), g_0'(s)) \, ds$$
$$= (\Phi^{x,a})^*(g_0, g_0).$$

This shows that condition 11 of Theorem 2.2 holds. Applying this result to $\{\mathcal{L}(Y_n^{x,a})\}_{n \in N}$, we have: for every set $A \in \mathcal{C}$,

(3.37) $$\liminf_n n^{-1} \log \mathbf{P}\{Y_n^{x,a} \in A\} \geq - \inf_{f \in A^0} (\Phi^{x,a})^*(f,f).$$

Finally, by a well-known argument (see [5], pages 518 and 519) based on (3.37), Lemma 3.7 and the fact that, for all $f \in E$,

$$\lim_{a \downarrow 0} (\Phi^{x,a})^*(f,f) = (\Phi^x)^*(f,f) = I^x(f),$$

we have

$$\liminf_n n^{-1} \log \mathbf{P}\{Y_n^x \in A\} \geq - \inf_{f \in A^0} I^x(f).$$

This completes the proof of the lower bound, and hence the proof of Theorem 3.1.

$\square$

COROLLARY 3.9. *Assume* (3.3) *and* (3.4). *Furthermore, assume that the initial value problem in* $\mathbf{R}^d$,

$$f(0) = x, \qquad f'(t) = \triangledown G(f(t), 0), \qquad t \in T,$$

*has a unique solution* $f_x$. *Then* $\{Y_n^x\}_{n \in N}$ *converges in probability to* $f_x$ *exponentially fast: for every* $\varepsilon > 0$, *there exists* $b > 0$ *such that*

$$\lim_n e^{bn} \mathbf{P}\{\|Y_n^x - f_x\|_\infty \geq \varepsilon\} = 0.$$

REMARK 3.10. As is well known (see, e.g., [12], page 270), a sufficient condition for the existence and uniqueness of $f_x$ is that the function $H(y) = \triangledown G(y, 0)$ satisfy a global Lipschitz condition on $\mathbf{R}^d$. This is closely related to condition (3.5).

PROOF OF COROLLARY 3.9. We claim first that

(3.38) $\qquad I^x(f) = 0 \qquad$ if and only if $f = f_x$.

For, it is easily seen that $G^*(y, z) = 0$ if and only if $z = \triangledown G(y, 0)$. The fact that $I^x(f_x) = 0$ is then clear. On the other hand, if

$$I_x(f) = \int_T G^*(f(t), f'(t)) \, dt = 0,$$

then $G^*(f(t), f'(t)) = 0$ a.e. $[m]$, and therefore $f'(t) = \triangledown G(f(t), 0)$ a.e. $[m]$. This implies that, for all $t \in T$,

$$f(t) = x + \int_0^t \triangledown G(f(s), 0) \, ds,$$

and therefore $f = f_x$. This proves (3.38).

Let $\varepsilon > 0$. By the upper bound statement of Theorem 3.1,

$$\limsup_n n^{-1} \log \mathbf{P}\{\|Y_n^x - f_x\|_\infty \geq \varepsilon\} \leq -\ell(\varepsilon),$$

where $\ell(\varepsilon) = \inf\{I^x(f) : f \in (B(f_x,\varepsilon))^c\}$. By the compactness of the level sets of $I^x$, there exists $g \in (B(f_x,\varepsilon))^c$ such that $I^x(g) = \ell(\varepsilon)$. Since $g \neq f_x$, (3.38) implies $I^x(g) > 0$. This establishes the conclusion. $\square$

Let $\mu$ be a probability measure on $\mathbf{R}^d$ such that $\hat{\mu}(\alpha) < \infty$ for all $\alpha \in \mathbf{R}^d$, where

$$\hat{\mu}(\alpha) = \int e^{\langle y,\alpha \rangle} \mu(dy).$$

Let $\{Z_k\}_{k \in N}$ be a sequence of i.i.d. random vectors in $\mathbf{R}^d$ with $\mathcal{L}(Z_1) = \mu$. Let $b : \mathbf{R}^d \to \mathbf{R}^d$ and $\sigma : \mathbf{R}^d \to \mathbf{R}^{d \times d}$ be bounded and uniformly Lipschitz. An interesting class of cases of the recursive scheme of Theorem 3.1 is obtained by taking

$$F_k(y) = b(y) + \sigma(y) Z_k, \qquad k \in \mathbf{N}.$$

Then we have, for $n \in \mathbf{N}, 0 \leq k \leq n$,

$$X^x_{n,0} = x,$$
$$X^x_{n,k} = X^x_{n,k-1} + n^{-1} b(X^x_{n,k-1}) + n^{-1} \sigma(X^x_{n,k-1}) Z_k, \qquad k \geq 1,$$

and $\{Y^x_n\}_{n \in N}$ may be regarded as a stochastic Euler polygonal scheme for the dynamical system

(3.39) $$f'(t) = b(f(t))$$

with initial condition $f(0) = x$.

THEOREM 3.11. (i) *If $\mu, b$ and $\sigma$ are as above, then $\{\mathcal{L}(Y^x_n)\}_{n \in N}$ satisfies the large deviation principle on $C(T, \mathbf{R}^d)$ with rate function*

$$I^x(f) = \begin{cases} \int_T G^*(f(t), f'(t)) \, dt, & \text{if } f(0) = x \\ & \text{and } f \text{ is absolutely continuous,} \\ \infty, & \text{otherwise,} \end{cases}$$

*where*

$$G^*(y, z) = \sup_{\alpha \in R^d} [\langle z, \alpha \rangle - G(y, \alpha)]$$

*and*

(3.40) $$G(y, \alpha) = \langle b(y), \alpha \rangle + \log \hat{\mu}(\sigma^t(y) \alpha),$$

*$\sigma^t(y)$ being the transpose of $\sigma(y)$.*

*Moreover, the level sets $L(I^x, h)$ are compact.*

(ii) $\{Y_n^x\}$ *converges in probability to the unique solution of the initial value problem*

$$f'(t) = b(f(t)) + \sigma(f(t))z_1, \qquad f(0) = x,$$

*exponentially fast* (*in the sense of Corollary* 3.9), *where* $z_1 = \mathbf{E}(Z_1)$.

REMARK 3.12. It is easily seen that if $\sigma(y)$ is invertible for all $y \in \mathbf{R}^d$, then

$$G^*(y, z) = (\log \hat{\mu})^*(\sigma^{-1}(y)(z - b(y))),$$

and therefore

$$I^x(f) = \begin{cases} \int_T (\log \hat{\mu})^*(\sigma^{-1}(f(t)))(f'(t) - b(f(t)))) \, dt, \\ \qquad \text{if } f(0) = x, \text{and } f \text{ is absolutely continuous,} \\ \infty, \qquad \text{otherwise.} \end{cases}$$

PROOF OF THEOREM 3.11. (i) We apply Theorem 3.1. Since in the present case

$$\int \exp(\langle z, \alpha \rangle) \mu(y, dz) = \mathbf{E} \exp\langle F_1(y), \alpha \rangle$$
$$= \mathbf{E} \exp[\langle b(y), \alpha \rangle + \langle \sigma(y)Z_1, \alpha \rangle]$$
$$= \exp[\langle b(y), \alpha \rangle] \hat{\mu}(\sigma^t(y)\alpha),$$

$G$ is indeed given by (3.40). Conditions (3.3) and (3.4) of Theorem 3.1 clearly hold. Since

$$(3.41) \quad \nabla G(y, \alpha) = b(y) + (\hat{\mu}(\sigma^t(y)\alpha))^{-1} \int \sigma(y)z \exp(\langle z, \sigma^t(y)\alpha \rangle) \mu(dz),$$

condition (3.5) follows from the assumptions on $\mu, b$ and $\sigma$. Finally,

$$|\langle F_1(y) - F_1(z), \alpha \rangle| = |\langle b(y) - b(z), \alpha \rangle + \langle (\sigma(y) - \sigma(z))Z_1, \alpha \rangle|$$
$$\leq C\|y - z\|\|\alpha\|(1 + \|Z_1\|)$$

for some constant $C > 0$ and therefore, for any $y, z \in \mathbf{R}^d, y \neq z, \tau > 0$,

$$\mathbf{E} \exp[\tau(\|y - z\|)^{-1}\langle F_1(y) - F_1(z), \alpha \rangle] \leq \mathbf{E} \exp(\tau C \|\alpha\|(1 + \|Z_1\|)).$$

This shows that condition (3.6) holds. By Theorem 3.1, statement (i) holds.

(ii) By (3.41),

$$\nabla G(y, 0) = b(y) + \int \sigma(y) z \mu(dz)$$
$$= b(y) + \sigma(y)z_1.$$

The statement follows now from Corollary 3.9. □

REMARK 3.13. 1. The particular case of (i) of Theorem 3.11 when $d = 1, \sigma(y) \neq 0$ for all $y \in \mathbf{R}$ and $\mu$ satisfies the additional assumptions $z_1 = 0$ and

$$\lim_{|\alpha| \to \infty} |\alpha|^{-1} \log \hat{\mu}(\alpha) = \infty$$

is presented in Theorem 2.1 of [10] by methods different from ours (actually, in [10] a more general dependence scheme is considered). However, the proof is incomplete: a convergence property of the rate function is used without justification ([10], pages 65 and 66).

2. The particular case of (i) of Theorem 3.11 when $\sigma(y)$ is invertible for all $y \in \mathbf{R}^d$ is implicitly covered by the presentation in [7], Proposition 6.3.4.

Department of Mathematics
Case Western Reserve University
Cleveland, Ohio 44106
USA
e-mail: add3@po.cwru.edu